\documentclass[11pt]{article}
\pdfoutput=1 
\usepackage[utf8x]{inputenc}
\usepackage[english]{babel}
\usepackage{graphicx}           
\usepackage{url}
\usepackage{eurosym}
\usepackage[plain]{algorithm}
\usepackage{algorithmic} 
\usepackage{color} 
\usepackage[margin=3.3cm]{geometry}
\usepackage{setspace}
\usepackage{datetime}
\usepackage{pifont}
\usepackage{tikz}
\usepackage{natbib}
\usepackage{hyperref}
\usepackage{amsmath, amsthm,amsfonts,amssymb,latexsym,stmaryrd}

\newcommand{\D}{\mathcal{D}}	

\newcommand{\EE}{\mathbb{E}}			
\newcommand{\PP}{\mathbb{P}}			
\newcommand{\RR}{\mathbb{R}}			
\newcommand{\NN}{\mathbb{N}}			

\newcommand{\taudiv}{b}							
\newcommand{\mmax}{M}    						
\newcommand{\mdiv}{m_{\textrm{\tiny\rm div}}}	

\newcommand{\MP}{\mathcal{N}}			
\newcommand{\op}{\mathcal{G}}			
\newcommand{\densiteIDEr}{m}				

\newcommand{\eqdef}     {\stackrel{{\textrm{\rm\tiny def}}}{=}}	
\newcommand{\umoins}{{\smash{u^{\raisebox{-1pt}{\scriptsize\scalebox{0.5}{$-$}}}}}}

\newcommand{\dif}{\mathrm{d}}
\newcommand{\rmd}{{{\textrm{\upshape d}}}}

\renewcommand{\tilde}{\widetilde}

\newcommand{\ind}{{{\mathrm\mathbf1}}}

\newtheorem{theorem}      {Theorem}[section]
\newtheorem{theorem*}     {theorem}
\newtheorem{proposition}  [theorem]{Proposition}

\newtheorem{lemma}        [theorem]{Lemma}

\newtheorem{examples}      [theorem]{Examples}
\newtheorem{remark}       [theorem]{Remark}

\newtheorem{corollary}    [theorem]{Corollary}

\newtheorem{hypothesis}   [theorem]{Assumption}
\newtheorem{hypotheses}   [theorem]{Assumptions}

\usepackage{ifthen}    
\newboolean{showComments}        
\setboolean{showComments}{true}  
\newcommand{\fnote}[1]
    {{\mbox{}\\\noindent\color{red}
    \rule{1cm}{2mm}\hfill  #1 \hfill\rule{1cm}{2mm}}
    \typeout{---------- #1 ------------}}

\newcommand{\dontforget}[1]
      {\ifthenelse {\boolean{showComments}} {{\color{red}(#1)}} {}}
\newcommand{\Dontforget}[1]
      {\ifthenelse {\boolean{showComments}} {\fnote{#1}} {}}

\begin{document}

\title{On the variations of the principal eigenvalue with respect to a parameter in growth-fragmentation models}

\author{Fabien Campillo$^{1,2}$ \and Nicolas Champagnat$^{3,4,5}$ \and Coralie Fritsch$^{6,3,4,5}$}

\footnotetext[1]{Inria, MATHNEURO, Montpellier, F-34095, France}

\footnotetext[2]{Institut Montpelli\'erain Alexander Grothendieck, Montpellier, F-34095, France}

\footnotetext[3]{Universit\'e de Lorraine, Institut Elie Cartan de Lorraine,
    UMR 7502, Vand\oe uvre-l\`es-Nancy, F-54506, France}

\footnotetext[4]{CNRS, Institut Elie Cartan de Lorraine, UMR
    7502, Vand\oe uvre-l\`es-Nancy, F-54506, France}

\footnotetext[5]{Inria, TOSCA, Villers-l\`es-Nancy, F-54600, France \protect \\ 
				E-mail: {fabien.campillo@inria.fr}, {nicolas.champagnat@inria.fr},
				 {coralie.fritsch@inria.fr}}
\footnotetext[6]{CMAP, \'Ecole Polytechnique, UMR CNRS 7641, route de Saclay, 91128 Palaiseau Cedex, France}

\maketitle

\begin{abstract}
We study the variations of the principal eigenvalue associated to a growth-frag\-men\-ta\-tion-death equation with respect to a parameter acting on growth and fragmentation. To this aim, we use the probabilistic individual-based interpretation of the model. We study the variations of the survival probability of the stochastic model, using a generation by generation approach. Then, making use of the link between the survival probability and the principal eigenvalue established in a previous work, we deduce the variations of the eigenvalue with respect to the parameter of the model.


\paragraph{Mathematics Subject Classification (MSC2010):}
35Q92, 45C05, 60J80, 60J85, 60J25, 92D25.
\end{abstract}


\section{Introduction}
\label{sec:introduction}

In biology, microbiology and medicine, diverse models are used to describe structured populations. For example the growth of a bacterial population or of tumor cells can be represented, in a constant environment, by the following growth-fragmentation-death equation \citep{doumic2007a, Doumic2010a, laurencot2009a, fredrickson1967a, sinko1967a, bell1967a, metz1986a}
\begin{align*} 
	\frac{\partial}{\partial t} \densiteIDEr_t(x)
	+\frac{\partial}{\partial x} \bigl(g(x)\,\densiteIDEr_t(x)\bigr)
	+ \bigl(\taudiv(x)+D \bigr)\,\densiteIDEr_t(x)
	=  
	2\,\int_x^\mmax 
		\frac{\taudiv(z)}{z}\, 
		q\left(z,\frac{x}{z}\right)\, 
		\densiteIDEr_t(z)\,\dif z\,,
\end{align*}
which describes the time evolution of the mass density $m_t$ of the population of cells which is subject to growth at speed $g$, cell division at rate $b$, with daughter cells generated by a division kernel $q$ and death at rate $D$. In order to study the asymptotic growth of the population, the eigenproblem associated to this equation is generally considered. The eigenvalue, also called Malthus parameter in this context, gives the asymptotic global growth rate of the population and allows to determine if the environment favors the development of the population.

Biologically, it is interesting to study the variation of this growth rate when its environment is changed (either by the action of an experimentalist or due to fluctuations of external conditions). In this article, we consider the model described previously, in which the growth function and the division rate depend on an environmental parameter $S$ describing the constant environment. 
The death rate is assumed independent of $S$ since we have in mind chemostat in which death is due to dilution at fixed rate.
This parameter can, for example, represent an external resource or the influence of other populations supposed to be at equilibrium. 
The study of the influence of this parameter on the growth of the population is a question of biological interest for a better understanding of the model, but also of numerical interest, for example, for the study of mutant invasions in adaptive dynamics problems \citep{fritsch2016a}.

This new question seems to be difficult to approach with standard deterministic mathematical tools where, up to our knowledge, no result is available except a study of the influence of asymmetric division by \cite{michel2006a, Michel2005a} and an asymptotical study of the influence of the parameters by \cite{calvez2012a}. See also the work of \cite{olivier2016a} for a study of the impact of the variability in cells' aging and growth rates as well as the one of \cite{clairambault2006a} for comparison of Perron eigenvalue (for constant in time birth and death rates) and Floquet eigenvalue (for periodic birth and death rates). The approach that we propose in this article uses the probabilistic interpretation of the growth-fragmentation-death equation under the form of a discrete stochastic individual-based model. 
This class of piecewise deterministic Markov processes is studied a lot, with a particular recent interest to the estimation of the parameters of the model \citep{doumic2015a,hoang2015a,hoffmann2015a}.
In this individual-based model, the growth of the population is determined by its growth rate, but also by its survival probability in some constant environment. The link between the eigenvalue of the deterministic model and the survival probability of the stochastic model, which correspond to two different definitions of the biological concept of invasion fitness \citep{Metz1992a, metz2008a}, was established by \cite{campillo2016a}. 
Our goal is to use this link to deduce variation properties of the eigenvalue with respect to the environmental parameter $S$ from the variations on the survival probability.
The probabilistic invasion fitness allows to use a generation by generation approach, which is more difficult to apply to the eigenproblem since generations overlap. Using this approach, the variations of the survival probability can be obtained by applying a coupling technique to the random process.

In an adaptive dynamics context, the variation of both invasion fitnesses are numerically very useful. For instance, considering the time evolution of a bacterial population in a chemostat, the invasion fitness determines if some mutant population can invade a resident one when a mutation occurs \citep{Metz1996a}. This invasion fitness is the one of the mutant population in the environment at the equilibrium determined by the resident one. In this example, the environmental parameter $S$ represents the substrate concentration at the equilibrium of the resident population. 
When the mutant population appears in the chemostat it appears in small size, hence its influence on the resident population and on the resource concentration can be neglected, which allows to assume the substrate concentration $S$ to be constant as long as the mutant population is small.
Moreover, due to the small number of mutant individuals, it is essential to use a stochastic model \citep{fritsch2015a, campillo2015b}.
However, the stochastic invasion fitness is numerically less straightforward to compute than the deterministic one. The mutual variations of both invasion fitnesses established in this article allow to considerably simplify the numerical analysis of a mutant invasion since the problem is reduced to the computation of a single eigenvalue in order to characterize the possibility of invasion of the mutant population \citep{fritsch2016a}.

In Section \ref{modele}, we present the deterministic and the stochastic versions of our growth-fragmentation-death model. We give the definitions of invasion fitness in both cases : for the stochastic one, it is defined as the survival probability and for the deterministic one, it corresponds to the eigenvalue of an eigenproblem. We extend some results from \cite{campillo2016a}, in particular Theorem \ref{th.critere} linking these two invasion fitnesses, to our more general context.
Section \ref{subsec.mon.sto} is devoted to the monotonicity properties of the survival probability of the stochastic model with respect to the initial mass and the death rate. In Section \ref{subsec.mon.sto.S} we prove, under suitable assumptions, the monotonicity of the survival probability with respect to the environmental parameter $S$. 
In Section \ref{subsec.mon.eigen}, we deduce from the previous results and from the link between the two invasion fitnesses, the monotonicity of the eigenvalue with respect to $S$.
Our assumptions are based on the realistic biological idea that the larger a bacterium is, the faster it divides and the larger the parameter $S$ is, the faster a bacterium grows. This is biologically consistent in the case where $S$ represents the substrate concentration. The monotonicity of fitnesses is obtained under additional assumptions which are detailed in the following sections. 
 We extend this result assuming a particular form of the growth rate $g$ and give a more general approach in Section \ref{subsec.ext.ccl}.

\section{Models description}
\label{modele}

In this Section we present two descriptions of the growth-fragmentation-death model. 
This model is the one studied by \cite{campillo2016a}, in which we add a dependence in a one-dimensional environmental parameter $S$, which is supposed to be fixed in time. 
In Section \ref{variation.environnement}, we study the variation of the invasion possibility of the population (whose definition depends on the considered description) with respect to $S$ for both descriptions.

\subsection{Basic mechanisms}
\label{subsec.mecha}

We consider models in which each individual is characterized by its mass $x \in [0,\mmax]$, where $\mmax$ is the maximal mass of individuals, and is affected by the following mechanisms:

\begin{enumerate}
\item \label{item.division} \textbf{Division:} each individual of mass $x$ divides at rate $\taudiv(S,x)$, into two individuals with masses $\alpha \, x$ and $(1-\alpha)\,x$, where the proportion $\alpha$ is distributed according to the probability distribution $Q(x,\dif \alpha)=q(x,\alpha)\,\dif \alpha$ on $[0,1]$.
\begin{center}
\includegraphics[width=5cm]{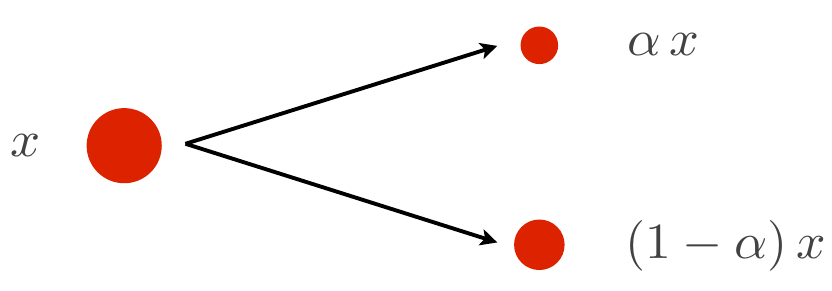}
\end{center}

\item \label{item.soutirage} \textbf{Death:} each individual dies at rate $D$.

\item \label{item.croissance} \textbf{Growth:} between division and death times, the mass of an individual grows at speed $g:\RR_+ \times [0,\mmax]\to\RR_+$ depending on an environmental parameter $S$, i.e.
\begin{align}
\label{eq.croissance}
  \frac{\rmd}{\rmd t} x_t
  = 
  g(S, x_t)\,.
\end{align}
\end{enumerate}

In this model, individuals do not interact between themselves and the environmental parameter $S$ is fixed in time. This means that the resource $S$ is not limiting for the growth of the population, this is for example the case if the resource is continuously kept at the same level or the consumption of the resource is negligible with respect to the resource quantity.
This model is relevant for a population with few individuals in a given environment such that the resource consumption is low.

\bigskip

For any $S>0$, let $A^S_t$ be the flow associated to an individual's mass growth in the environment $S$, i.e. for any $x\in (0,M)$ and $t\geq 0$,
\begin{align}
\label{def.flot}
	A^S_t(x) = x + \int_0^t g(S, A^S_u(x))\,\dif u \,.
\end{align}

Throughout this paper we assume the following set of assumptions.

\begin{hypotheses}
\label{da.hypo.model.reduit}
\begin{enumerate}
\item \label{hyp.q}For any $x\in[0,\mmax]$, the kernel $q(x,.)$ is symmetric with respect to $1/2$:
$$
	q(x,\alpha) = q(x,1-\alpha), \qquad \forall \alpha \in [0,1]
$$
such that $\int_0^1 q(x,\alpha)\,\dif \alpha = 1$.
\item \label{hyp.q.continue} For any $\alpha \in [0,1]$, the function $x\mapsto q(x,\alpha)$ is continuous on $[0,\mmax]$.

\item \label{hyp.q.domine} There exists a function $\bar q : [0,1] \mapsto \RR_+$ such that $q(x,\alpha) \leq \bar q (\alpha)$ for any $x\in(0,\mmax)$ and $\int_0^1 \bar q(\alpha)\,\dif \alpha < + \infty$.

\item \label{g.0.M} $g(S,0)=g(S,M)=0$ and  $g(S,x)>0$ for any $x \in (0,M)$ and $S>0$.

\item \label{hyp.g} $g(S,.)\in C[0,\mmax]\cap C^1(0,\mmax)$ , where $C[0,\mmax]$ and $C^1(0,\mmax)$ respectively represent sets of continuous functions on $[0,\mmax]$ and continuously differentiable functions on $(0,\mmax)$.

\item $b(S,.)\in C[0,\mmax]$ and there exists $\mdiv \in [0,M)$ and $\bar \taudiv>0$ such that
\begin{align*}
	\taudiv(S,x)=0 & \textrm{ if } x \leq \mdiv \,, 
\\
	0<\taudiv(S,x)\leq\bar \taudiv & \textrm{ if } x \in (\mdiv,M)\,.
\end{align*}
\end{enumerate}
\end{hypotheses}

Assumptions \ref{da.hypo.model.reduit}-\ref{hyp.g} and \ref{da.hypo.model.reduit}-\ref{g.0.M}
ensure existence and uniqueness of the growth flow defined by \eqref{def.flot} for $x\in(0,\mmax)$ until the exit time $T_{\text{exit}}(x):=\inf\{t>0 \, | \, A^S_t(x) \geq \mmax\}$ of $(0,\mmax)$ and that $A^S\in C^1(\D)$ with $\D=\{(t,x),\, t<T_{\text{exit}}(x) \}$ \cite[Th. 6.8.1]{demazure2000a}.
We define this flow as constant when it starts from $\mmax$.
Note that the exit time $T_{\text{exit}}(x)$ is infinite if the convergence $\lim_{x\to \mmax}g(S,x)=0$ is sufficiently fast (see for example \cite[Assumption 3.]{campillo2016a} for more details).
Assumption \ref{g.0.M} means that the maximal biomass of an individual is the same for any concentration of resources. This may not
  be true in general, but we can always change the scale of biomass for each value of $S$ so that the maximal value of $x$ is always
  $M$ and modify the growth and birth parameters accordingly. This is what we shall assume in the sequel.

\subsection{Growth-fragmentation-death integro-differential model}
\label{subsec.edp}

The deterministic model associated to the previous mechanisms is given by the integro-differential equation
\begin{align} 
\label{da.eid.reduit}
	\frac{\partial}{\partial t} \densiteIDEr^S_t(x)
	+\frac{\partial}{\partial x} \bigl(g(S, x)\,\densiteIDEr^S_t(x)\bigr)
	+ \bigl(\taudiv(S,x)+D \bigr)\,\densiteIDEr^S_t(x)
	=  
	2\,\int_x^\mmax 
		\frac{\taudiv(S,z)}{z}\, 
		q\left(z,\frac{x}{z}\right)\, 
		\densiteIDEr^S_t(z)\,\dif z\,,
\end{align}
where $\densiteIDEr^S_t(x)$ represents the density of individuals with mass $x$ at time $t$ evolving in the environment determined by $S$, with a given initial condition $\densiteIDEr^S_0$. 

\medskip

Let $\op_S$ be the non local transport operator such that $\partial_t m^S_t(x) = \op_S m^S_t(x)$: for any $f\in C^1(0,\mmax)$, $x\in(0,\mmax)$,
\begin{align}
\label{def.G}
  \op_S f(x)
  \eqdef
  - \partial_x (g(S, x)\,f(x))
  - (D+\taudiv(S,x))\, f(x)
  + 2\int_x^\mmax \frac{\taudiv(S,z)}{z}\, 
		q\left(z,\frac{x}{z}\right)\, 
		f(z)\,\dif z\,,
\end{align}
and $\op_S^*$ its adjoint operator defined for any $f\in C^1(0,\mmax)$, $x\in(0,\mmax)$ by
\begin{align}
\label{def.G.star}
  \op_S^* f(x)
  \eqdef
  - (D+\taudiv(S,x))\, f(x)
  + g(S,x)\,\partial_x f(x)
  + 2\,\taudiv(S,x)\,\int_0^1q(x,\alpha)\,f(\alpha\, x)\,\dif \alpha\,.
\end{align}

We consider the eigenproblem
\begin{subequations}
\label{eq.eigenproblem}
\begin{align}
\op_S \hat u_S(x) = \Lambda_S \, \hat u_S(x)\,,
\end{align}
\begin{align}
   \lim_{x\to 0}g(S,x)\, \hat u_S(x) &= 0\,,
   & 
   D+\Lambda_S &>0\,,
   &
   \hat u_S(x) &\geq 0\,,
   & 
   \int_0^\mmax \hat u_S(x)\,\dif x &= 1
\end{align}
\end{subequations}
and the adjoint problem
\begin{align}
\label{eq.eigenproblem.dual}
\op_S^* \hat v_S(x) = \Lambda_S \, \hat v_S(x)\,,
 \qquad
\hat v_S(x)\geq 0\,,
\qquad
   \int_0^\mmax \hat v_S(x)\,\hat u_S(x)\,\dif x = 1\,.
\end{align}

\medskip

The eigenvalue $\Lambda_S$ is then interpreted as the exponential growth rate (or decay rate if it is negative) of the population.

\medskip

In the rest of the paper, we will assume that the following assumption is satisfied. \cite{campillo2016a} have given some conditions under which this assumption holds (see also \citep{doumic2007a, Doumic2010a} for sligthly different models and \citep{perthame2005a, laurencot2009a, mischler2016a} for exponential stability of the eigenfunctions).

\begin{hypothesis}
\label{da.hyp.solution.pp}
For any $S>0$, the system \eqref{eq.eigenproblem}-\eqref{eq.eigenproblem.dual} admits a solution $(\hat u_S,\hat v_S, \Lambda_S)$ such that $\hat u_S \in C^1(0,\mmax)$ and $\hat v_S \in C[0,\mmax] \cap C^1(0,\mmax)$. 
\end{hypothesis}

\subsection{Growth-fragmentation-death individual-based model}
\label{subsec.IBM}

The mechanisms described in Section \ref{subsec.mecha} can also be represented by a stochastic individual-based model, where the population at time $t$ is represented by the counting measure
\begin{align}
\label{da.pop.modele.reduit}
	\eta^S_t(\rmd x) \eqdef \sum_{i=1}^{N_t}\delta_{X_t^i}(\rmd x)\,,
\end{align}
where $N_t=\int_0^\mmax \eta_t^S(\dif x)$ is the number of individuals in the population at time $t$ and $(X_t^i,\, i=1,\dots, N_t)$ are the masses of the $N_t$ individuals (arbitrarily ordered).

The stochastic individual-based model is relevant for small population whereas the deterministic one is relevant for large population \citep{campillo2015b}.

The process $(\eta^S_t)_{t\geq 0}$ is defined by
\begin{align}
\nonumber
  &\eta^S_{t}
  =
  \sum_{j=1}^{N_0}\delta_{A^S_t(X_0^j)}
\\
\nonumber
  &
  \qquad+
  \iiiint\limits_{[0,t]\times\NN^*\times[0,1]^3}
        1_{\{j\leq N_{\umoins}\}} \, 
        1_{\{\theta_1 \leq \taudiv(S,X_\umoins^j)/\bar \taudiv\}}\,  
        1_{\{\theta_2 \leq q(X_\umoins^j,\alpha)/\bar q(\alpha)\}}\,
        \bigl[
                -\delta_{A^S_{t-u}(X^j_\umoins)}
\\[-1.2em]
\nonumber
&\qquad\qquad\qquad\qquad\quad
                +\delta_{A^S_{t-u}(\alpha\,X^j_\umoins)}	
                +\delta_{A^S_{t-u}((1-\alpha)\,X^j_\umoins)}
         \bigr]\,
		 		\MP_1(\dif u, \dif j, \dif \alpha, \dif \theta_1, \dif \theta_2) 
\\[0.2em]
\label{da.def.proc.S.nu}
  &\qquad-
		\iint\limits_{[0,t]\times\NN^*} 
			1_{\{j\leq N_{\umoins}\}} \, 
			 \delta_{A^S_{t-u}(X^j_\umoins)}
			 			 \,	\MP_2(\dif u, \dif j)\, 
\end{align}
where $\MP_1(\dif u, \dif j, \dif \alpha, \dif \theta_1,\theta_2)$ and $\MP_2(\dif u, \dif j)$ are two independent Poisson random measures defined on $\RR_+ \times \NN^* \times [0,1] \times [0,1]\times [0,1]$ and $\RR_+ \times \NN^*$, corresponding respectively to the division and death mechanisms, with respective intensity measures
\begin{align}
n_1(\dif u, \dif j, \dif \alpha, \dif \theta)
	& =
		\bar \taudiv \, \dif u \, \Big(\sum_{\ell \geq 1} \delta_{\ell}(\dif j) \Big)
		\, \bar q(\alpha)\, \dif \alpha \, \dif \theta_1\, \dif \theta_2\,,
\\
n_2(\dif u, \dif j)
	& =
		D \, \dif u \, \Big(\sum_{\ell \geq 1} \delta_{\ell}(\dif j) \Big)\,,
\end{align}
(see \cite{campillo2015b} and \cite{campillo2016a} for more details).

This population process can be seen as a multitype branching process with a continuum of types.
We are interested in its survival probability.

We suppose that, at time $t=0$, there is only one individual, with mass $x_0$, in the population, i.e.
$$
	\eta^S_0(\rmd x) = \delta_{x_0}(\rmd x)\,.
$$
 
\bigskip

The extinction probability of the population with initial mass $x_0$ is
\begin{align*}
  p^S(x_0) \eqdef \PP^S_{\delta_{x_0}}(\exists t>0, N_t=0)\,,
\end{align*}
where $\PP^S_{\delta_{x_0}}$ is the law of the process $(\eta_t^S)_{t\geq 0}$ under the initial condition $\eta^S_0=\delta_{x_0}$.
The survival probability is then given by $\PP^S_{\delta_{x_0}}(\text{survival})=1- p^S(x_0)$.

We define the $n$-th generation as the set of individuals descended from a division of one individual of the $(n-1)$-th generation.
The generation 0 corresponds to the initial population.
We denote by $Z_n$ the number of individuals at the $n$-th generation and
we define the extinction probability before the $n$-th generation as
\begin{align*}
  p^S_n(x_0) 
  \eqdef
  \PP^S_{\delta_{x_0}}
    (Z_n=0)\,,
    \quad n\in \NN\,.
\end{align*}
It is obvious that
$$
	\lim_{n \to \infty} p^S_n(x_0)=p^S(x_0)\,.
$$

\medskip

Let $\tau$ be the stopping time of the first event (division or death). Then at time $\tau$ the population is given by
\begin{align}
\label{eta_tau}
  \eta^S_\tau
  \eqdef
  \left\{
			\begin{array}{ll}
				0\, & \text{if death}\,,\\
				\delta_{X_1} + \delta_{X_2}\, & \text{if division}\,,
			\end{array}
  \right.
\end{align}
with $X_1=\alpha \, A^S_\tau(x_0)$ and $X_2=(1-\alpha)\,A^S_\tau(x_0)$ where the proportion $\alpha$ is distributed according to the kernel $q(A^S_\tau(x_0),\alpha)\,\dif \alpha$.

\bigskip

Applying the Markov property at time $\tau$ and using the independence of particles, it is easy to prove (see \cite{campillo2016a}) that for any $x\in [0,\mmax]$ and $n\in \NN^*$
\begin{multline}
\label{proba.extinction.n.gen}
  p_n^S(x)
  =
  D\, \int_0^\infty e^{-D\,t}\, e^{-\int_0^t \taudiv(S,A^S_u(x))\,\dif u}\, 
  \dif t
  + \int_0^\infty \taudiv(S,A^S_t(x)) \, e^{-\int_0^t \taudiv(S,A^S_u(x)) 
      \,\dif u -D\,t}\,
\\
	\int_0^1 q(A^S_t(x),\alpha)\,p^S_{n-1}\big(\alpha \, A^S_t(x)\big)\,
	     p^S_{n-1}\big((1-\alpha) \, A^S_t(x)\big)\, \dif \alpha\, \dif t\,.
\end{multline}
with $p^S_0(x)=0$. It can then be deduced \cite[Proposition 3]{campillo2016a} that $p^S$ is the minimal non negative solution of
\begin{align}
\nonumber
  p^S(x)
  &=
  \int_0^\infty D\, e^{-D\,t}\, e^{-\int_0^t \taudiv(S,A^S_u(x))\,\dif u}\, 
  \dif t
\\
\nonumber
  &\quad
  + \int_0^\infty \taudiv(S,A^S_t(x)) \, e^{-\int_0^t \taudiv(S,A^S_u(x)) 
      \,\dif u -D\,t}\,
\\
\label{da.eq.proba.extinction}
  &\qquad\qquad \qquad
	\int_0^1 q(A^S_t(x),\alpha)\,p^S\big(\alpha \, A^S_t(x)\big)\,
	     p^S\big((1-\alpha) \, A^S_t(x)\big)\, \dif \alpha\, \dif t \,,
\end{align}
in the sense that for any non negative solution $\tilde p$ we have $\tilde p\geq p^S$.

\medskip

\begin{remark}
By a change of variable, we have 
\begin{align*}
  p^S(x)
  &=
  \int_x^\mmax \frac{D}{g(S,y)}\, e^{-\int_x^y \frac{\taudiv(S,z)+D}{g(S,z)}\,\dif z}\, 
  \dif y
\\
  &\quad
  + \int_x^\mmax \frac{\taudiv(S,y)}{g(S,y)} \, e^{-\int_x^y \frac{\taudiv(S,z)+D}{g(S,z)}\,\dif z} 
	\int_0^1 q(y,\alpha)\,p^S\big(\alpha \, y\big)\,
	     p^S\big((1-\alpha) \, y\big)\, \dif \alpha\, \dif y \,.
\end{align*}
Therefore, the extinction probability is solution of
$$
g(S,x)\partial_x p^S(x) + D\,(1-p^S(x))
+ b(S,x)\,\Big\{\int_0^1 q(x,\alpha)\,p^S(\alpha\,x)\,p^S((1-\alpha)\,x)\,\dif \alpha
-p^S(x)\Big\}
=0\,.
$$
\end{remark}

For any $x \in ]0,M[$ and $y>0$ such that $x\leq y$, let $t_S(x,y)$ be the first hitting time of $y$ by the flow $A^S_t(x)$, i.e. 
\begin{align}
\label{temps.atteinte.masse}
  t_S(x,y)
  \eqdef
  \inf\{t \geq 0 ,\, A^S_t(x)=y \}
  =
  \left\{
			\begin{array}{ll}
				\tilde A_{S,x}^{-1}(y) \,, & \text{if $x\leq y<M$}\,,\\
				+\infty\,, & \text{if $y\geq M$}\,,
			\end{array}
  \right.
\end{align}
where $\tilde A_{S,x}^{-1}$ is the inverse function of the $C^1$-diffeomorphism $t\mapsto A^S_t(x)$.

\medskip

\cite{campillo2016a} we have made the link between the survival probability of the stochastic process and the eigenvalue of the
deterministic model, given by the theorem below. This result was proved for a kernel $q(x,.)$ which does not depend on $x\in
(0,\mmax)$, but it can easily be extended to our case where $q(x,.)$ depends on the mass $x$ at the division time as explained below.

\medskip

\begin{theorem}[Campillo, Champagnat, Fritsch (2015)]
\label{th.critere}
Under Assumptions \ref{da.hypo.model.reduit} and \ref{da.hyp.solution.pp}, we have the following relation between the two invasion criteria
$$\Lambda_S>0 \quad \Longleftrightarrow \quad 
\PP^S_{\delta_{x}}(\text{survival})>0\,,\, \forall x\in(0,\mmax)\,.$$ 
\end{theorem}

Note that, contrary to the works of \cite{perthame2005a,doumic2007a,laurencot2009a,Doumic2010a,mischler2016a}, we assume here a compact set $[0,M]$ of biomasses to keep things simple in the sequel.
  The extension of our approach to a non-compact case would require to identify the good assumptions at infinity for the last result
  to hold (the rest of our arguments should work similarly). The last problem is not so easy because it strongly depends on the
  growth at infinity of the eigenfunctions $\hat{u}$ and $\hat{v}$ of Assumption~\ref{da.hyp.solution.pp}. Note in addition that the
  problem of existence of these eigenfunctions also requires a careful study at infinity (see~\cite{Doumic2010a}).

\medskip

\begin{proof}
The key argument of the proof is that the process $(e^{-\Lambda_S\,t}\,\sum_{i=1}^{N_t}\hat v_S(X_t^i))_{t\geq 0}$ is a $\PP^S_{\delta_x}$-martingale such that 
$$
e^{-\Lambda_s\,t}\,\sum_{i=1}^{N_t}\hat v_S(X_t^i)
\xrightarrow[t\to\infty]{} \mathcal{Z}\ \textrm{ $\PP^S_{\delta_{x}}$-a.s.}
$$
where $\mathcal{Z}$ is an integrable random variable (see \cite[Theorem 2 and Lemma 3]{campillo2016a}). The arguments of \cite{campillo2016a} to prove that
\begin{enumerate}
\item if $\Lambda_S>0$ then $\PP^S_{\delta_x}(\text{survival})>0$ for any $x\in(0,\mmax)$
\item if $\Lambda_S<0$ then $\PP^S_{\delta_x}(\text{survival})=0$ for any $x\in[0,\mmax]$
\end{enumerate}
can be directly applied for a kernel $q(x,.)$ depending on the mother mass $x$. The first  statement is proved using that if $\Lambda_S>0$ then $\mathcal{Z}$ is bounded in $L^2$ whereas the second one comes from the inequality
$$
	\EE_{\delta_x}(N_t) \leq C_x\,e^{\Lambda_S\,t}, \qquad \forall t\geq 0\,
$$
where $C_x>0$ is a constant depending on the initial mass $x\in (0,\mmax)$. Its proof by~\cite{campillo2016a} is
  technical, but the extension to kernels $q$ depending on $x$ is straightforward.

The only difficulty concerns the third point of the proof of \cite[Theorem 2]{campillo2016a} in which we prove that if $\Lambda_S =
0$ then $M_t^\varepsilon\to 0$ a.s. with $M_t^\varepsilon$ the number of individuals with mass in $[\varepsilon,\mmax-\varepsilon]$
at time $t$ for $0<\varepsilon<\frac{\mmax}{2}$. Then, the fourth point of the proof, stating that $M_t^\varepsilon\to 0$ implies
extinction, follows similarly. The main idea of this fourth step is that the number of individuals $M_t^\varepsilon$ cannot
indefinitely stay in a compact subset $\{1,2,\ldots,c\}$ of $\mathbb{N}$. Then either $\limsup_{t\to\infty} M_t^\varepsilon =
\infty$ or $M_t^\varepsilon \to 0$ a.s.\ However, $\limsup_{t\to\infty} M_t^\varepsilon = \infty$ contradicts the fact that
$\mathcal{Z}$ is integrable if $\Lambda_S=0$, hence $M_t^\varepsilon \to 0$ a.s.

For the proof of the third point (see details in \citep{campillo2016a}), it is sufficient that for $c>0$, there exists $t_0$ such that
\begin{align}
\label{eq.gamma}
1>
\gamma
	:=
		\inf_{\varepsilon \leq x \leq \mmax-\varepsilon} 
		\PP^S_{\delta_{x}}\left(M_{t_0}^\varepsilon \geq c \right)>0\,.
\end{align}
If $q(x,.)=q(.)$ is independent of the mother mass, it is sufficient to take $\varepsilon>0$ such that $x\in(\varepsilon,\mmax-\varepsilon)$ and 
$\textstyle q\left(\left[ \varepsilon / (\mmax-2\,\varepsilon), \, 1/2\right]\right)>0$ to obtain \eqref{eq.gamma} for one $t_0$. 
In our case, the last condition must be replaced by $\inf_{\varepsilon\leq x\leq
  \mmax-\varepsilon}q(x,[\varepsilon/(\mmax-2\,\varepsilon),1/2])>0$  for some $\varepsilon>0$. 
Note that the infimum above is reached at some $x_0(\varepsilon) \in [\varepsilon,\mmax-\varepsilon]$, by Assumptions \ref{da.hypo.model.reduit}-\ref{hyp.q.continue} and \ref{hyp.q.domine}.
Therefore, we proceed by contradiction and assume that for all $\varepsilon>0$, there exists $x_0(\varepsilon)\in[\varepsilon,M-\varepsilon]$  such that $q(x_0(\varepsilon),\alpha)=0$ for almost all $\alpha\in\left[\frac{\varepsilon}{\mmax-2\,\varepsilon} \,; \, \frac 12\right]$. Then, from the sequence $\left(x_0\left(\frac 1n\right)\right)_n$, we can extract a subsequence which converges towards $x_0^*$. By continuity of $x\mapsto q(x,\alpha)$, we then get $q(x_0^*,\alpha)=0$ for almost all $\alpha \in (0,1)$. Hence $\int_0^1 q(x_0^*,\alpha)\,\dif \alpha = 0$, which contradicts Assumption \ref{da.hypo.model.reduit}-\ref{hyp.q}.
\end{proof}

\section{Variations of the invasion fitnesses with respect to the environmental variable}
\label{variation.environnement}

Our goal is to study the variation of $\Lambda_S$ w.r.t $S$. For this, we start by studying the monotonicity properties of the survival probability in the stochastic model.

\subsection{Monotonicity properties w.r.t. the initial mass and the death rate on the stochastic model}
\label{subsec.mon.sto}

  From a biological point of view, little is known about the dependence of the division kernel $q$ w.r.t.\ $x$ \citep{osella2017a}. Most
  often, it is assumed independent of $x$ in applications. In order to obtain the most general result, we assume that $q$ depends on
  $x$, and we need to state assumptions about this parameter. Note however that the self-similar fragmentation is included in our assumptions. Moreover, although $q(x,\alpha)$ is assumed to be regular w.r.t. $\alpha$, more general kernels can be consider, in particular the following results should hold for self-similar equal mitosis.

For any $x\in(0,1)$, let $F_x : [0,1] \to [0,1]$ be the cumulative distribution function associated to the law $q(x,\alpha)\,\dif \alpha$, that is for any $u\in[0,1]$
$$
	F_x(u) = \int_0^u q(x,\alpha)\,\dif \alpha
$$
and let $F_x^{-1}$ be its inverse function defined by
$$
	F_x^{-1}(v) = \inf_{u\in [0,1]} \left\{F_x(u)\geq v \right\}\,.
$$

\begin{hypothesis}
\label{hyp.fct.rep.q}
The cumulative distribution function $F_x$ satisfies, for any $u\in(0,1)$ and any $x\leq y$,
$$
	x\,F_x^{-1}(u) \leq y \,F_y^{-1}(u)  \qquad \text{and} \qquad
	(1-x)\,F_x^{-1}(u) \leq (1-y) \,F_y^{-1}(u) \,.
$$
\end{hypothesis}

As we will see in Lemma \ref{lemma.int.q.dec} below, this assumption corresponds to a coupling condition on the mass of offspring
born from individuals of different sizes. We need this condition because our method can be seen as a construction of a coupling of
the masses of individuals at each generation in two stochastic processes starting from different initial masses (see our comments below, particularly Remark~\ref{rk:proba-prop-3.5}).

\medskip

\begin{remark}
\label{remark.cdf}
If $F_x^{-1}$ is such that for any $u\in(0,1)$, $x\mapsto F_x^{-1}(u) \in C^1([0,\mmax])$ and satisfies for any $x\in(0,\mmax)$,
$$
	x\,\partial_x F_x^{-1}(u) \in [-F_x^{-1}(u), 1-F_x^{-1}(u)] \,,
$$
then
$$
	\partial_x(x\,F_x^{-1}(u)) = F_x^{-1}(u) + x\,\partial_x F_x^{-1}(u) \geq 0\,.
$$
Hence $x\,F_x^{-1}(u)$ is non decreasing. In the same way, $(1-x)\,F_x^{-1}(u)$ is non decreasing too. Therefore, Assumption \ref{hyp.fct.rep.q} holds.
\end{remark}

\medskip

\begin{figure}
\begin{center}
\includegraphics[scale=0.24]{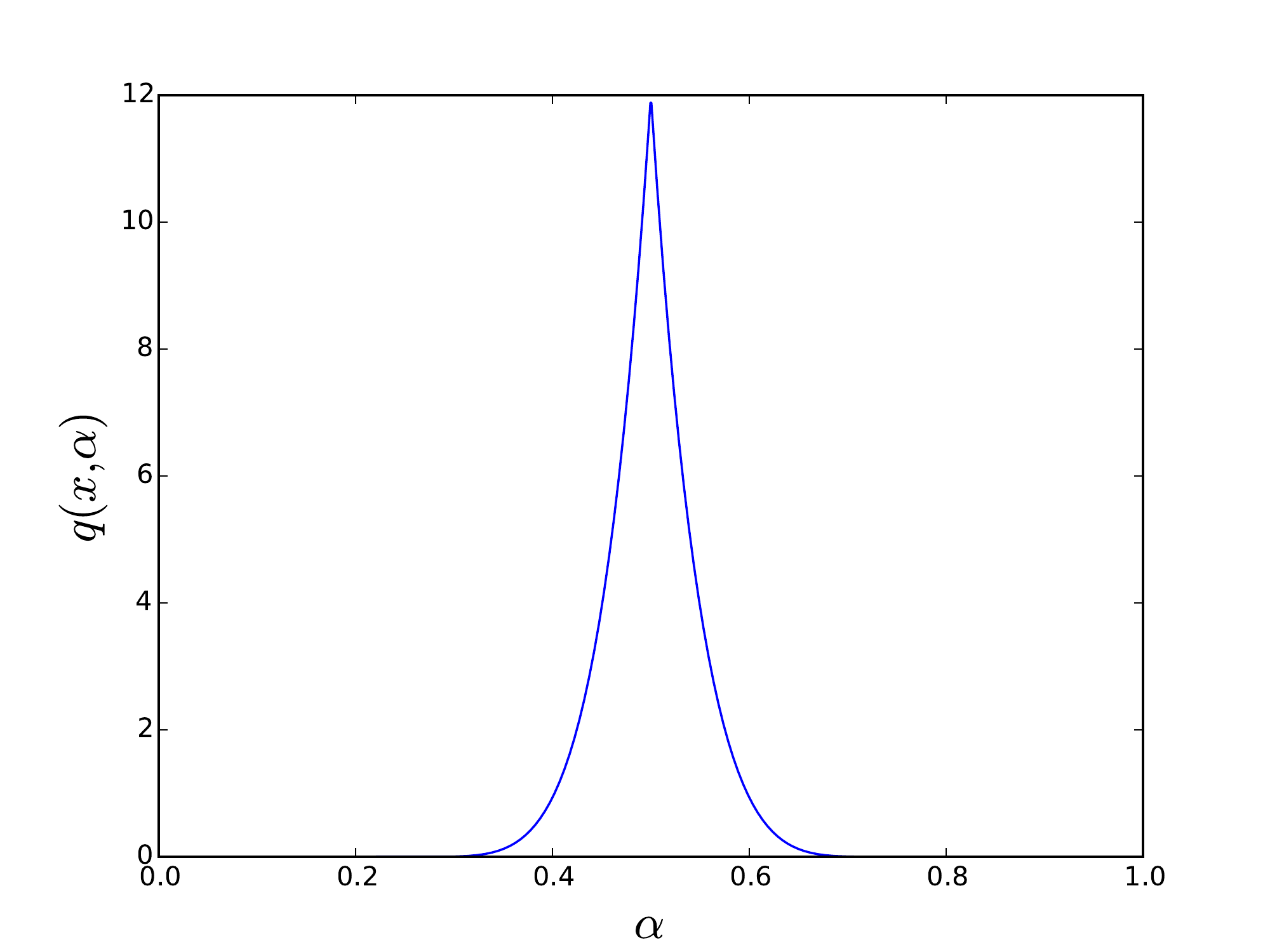}
\includegraphics[scale=0.24]{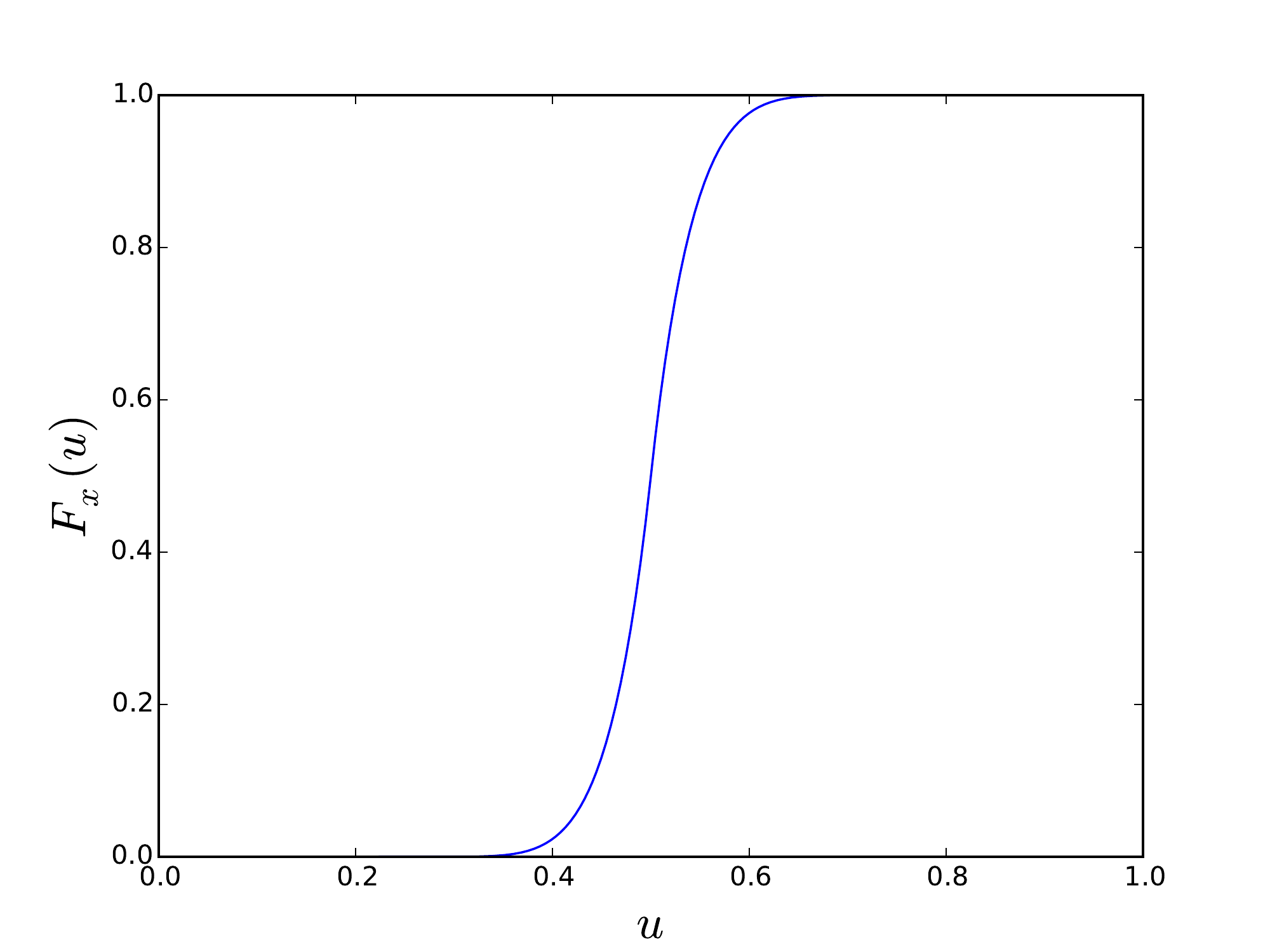}
\includegraphics[scale=0.24]{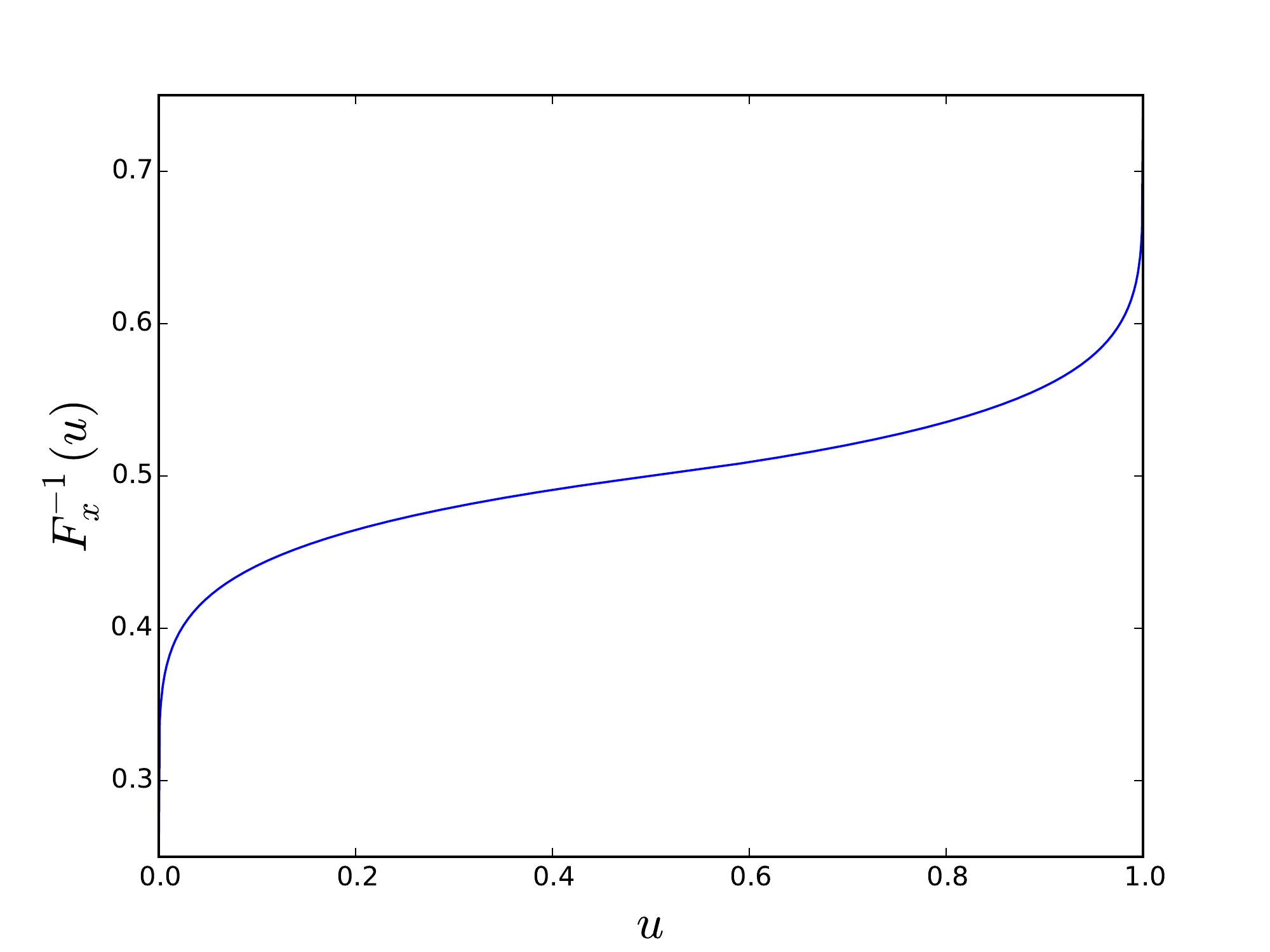}
\end{center}
\caption{\label{fig.ex.q} Representation of the function $q$ (left), $F_x$ (center) and $F_x^{-1}$ (right) respectively defined by Equations \eqref{eq.ex.q}, \eqref{eq.ex.F} and \eqref{eq.ex.F.inv} with $l(x)=0.25$ and $\beta(x)=5$.}
\end{figure}

\begin{examples}
We give some examples which satisfy Assumption \ref{hyp.fct.rep.q}.
\begin{enumerate}
\item We consider the following division kernel,
\begin{align*}
	q(x,\alpha) = \frac{\ind_{\{l(x) \leq \alpha \leq 1-l(x)\}}}{1-2\,l(x)}\,.
\end{align*}
where $l \in C^1([0,\mmax],(0,1/2))$.
Then for $u \in (0,1)$,
$$
	F_x^{-1}(u) = (1-2\,u)\,l(x)+u
$$
and, by Remarks \ref{remark.cdf}, Assumption \ref{hyp.fct.rep.q} holds if for any $x$, $0\leq x\,l'(x)+l(x)\leq 1$.
\item We can extend the previous example considering the following function $q$,
\begin{align}
\label{eq.ex.q}
	q(x,\alpha) = \frac{(\alpha-l(x))^{\beta(x)}}{C(x)} \, \ind_{\{l(x)\leq \alpha \leq 1/2\}}
				+ \frac{(1-\alpha-l(x))^{\beta(x)}}{C(x)} \, \ind_{\{1/2\leq \alpha \leq 1-l(x)\}}
\end{align}
where $C(x) = 2\,\left(1/2-l(x) \right)^{\beta(x)+1}/(\beta(x)+1)$ is a normalizing constant.
The previous example corresponds to $\beta(x)=0$ for any $x\in[0,\mmax]$.
Then 
\begin{align}
\nonumber
F_x(u) 
	&= 
		\frac 12 \,
		\left(\frac{u-l(x)}{\frac 12-l(x)}\right)^{\beta(x)+1}\, 
		\! \! \!\! \! \! \! \! \! \ind_{\{l(x)\leq \alpha \leq 1/2\}}
\\
\label{eq.ex.F}
	& \quad 
	+\left(1-\frac 12 \,
				\left(\frac{1-u-l(x)}{\frac 12-l(x)}\right)^{\beta(x)+1}\right)
			\, \ind_{\{1/2< \alpha \leq 1-l(x)\}}
	+ \ind_{\{1-l(x)< \alpha\}}
\end{align}
and for any $u\in (0,1)$
\begin{align}
\nonumber
F^{-1}_x(u)
		& =
		\left(\left(\frac 12-l(x)\right)\,(2\,u)^{1/(\beta(x)+1)}+l(x)\right)\,
		\ind_{\{0< u \leq 1/2\}}
\\
\label{eq.ex.F.inv}
	& \quad
		+\left(1-l(x)-
			\left(\frac 12-l(x)\right)\,
			(2\,(1-u))^{1/(\beta(x)+1)}\right)\,\ind_{\{1/2< u < 1\}}\,.
\end{align}
An example of such functions is given in Figure \ref{fig.ex.q}.

For $u \in (0,1/2]$,
\begin{align*}
\partial_x F_x^{-1}(u)
	&=
		\left(
			- l'(x) -\left(\frac 12 - l(x)\right)\,
			\frac{\beta'(x)}{(\beta(x)+1)^2}\,\ln(2\,u)
		\right)\,(2\,u)^{1/(\beta(x)+1)}
		+l'(x)
\end{align*}
and for $u \in [1/2,1)$,
\begin{align*}
\partial_x F_x^{-1}(u)
	&=
	\left(
		l'(x)+\left(\frac 12 - l(x)\right)\,\frac{\beta'(x)}{(\beta(x)+1)^2}\,\ln(2\,(1-u))
	\right)\,(2\,(1-u))^{1/(\beta(x)+1)}
	-l'(x)
\end{align*}
Assumption \ref{hyp.fct.rep.q} holds if $0 \leq x\,\partial_x F_x^{-1}(u)+F_x^{-1}(u) \leq 1$ for any $u\in (0,1)$, for example if $\beta$ is a constant function and if $0\leq l(x)+ x\,l'(x)\leq 1$ for any $x \in (0,\mmax)$.
\end{enumerate}
\end{examples}

\medskip

\begin{lemma}
\label{lemma.int.q.dec}
Let $f$ be a non-increasing function on $[0,\mmax]$. Then, under Assumption \ref{hyp.fct.rep.q}, the function
$$
	x \mapsto \int_0^1 q(x,\alpha)\,f(\alpha\,x)\,f((1-\alpha)\,x)\,\dif \alpha
$$
is non-increasing.
\end{lemma}

\medskip

\begin{proof}
For any $x\in(0,1)$, let $\theta_x$ defined by $\theta_x=F_x^{-1}(U)$ where $U$ is uniformly distributed on $[0,1]$. Therefore the law of the variable $\theta_x$ is $q(x,\alpha)\,\dif \alpha$.
By Assumption \ref{hyp.fct.rep.q},
$$
	\partial_x(x\,\theta_x) = \theta_x + x\,\partial_x F_x^{-1}(U)
	\geq \theta_x - F_x^{-1}(U)=0 \quad \text{a.s.}
$$
and
$$
	\partial_x(x\,(1-\theta_x)) = 1-\theta_x - x\,\partial_x F_x^{-1}(U)
	\geq 1-\theta_x -(1-F_x^{-1}(U))=0\,.
$$
Therefore, for any $x<y$ we have $x\,\theta_x \leq y\,\theta_y$ a.s. and 
$x\,(1-\theta_x) \leq y\,(1-\theta_y)$ a.s. Hence,
\begin{align*}
\int_0^1 q(x,\alpha)\,f(\alpha\,x)\,f((1-\alpha)\,x)\,\dif \alpha
	&=
		\EE\left(f(\theta_x\,x)\,f((1-\theta_x)\,x) \right)
\\
	&\leq 
		\EE\left(f(\theta_y\,y)\,f((1-\theta_y)\,y) \right)
\\
	& =
	\int_0^1 q(y,\alpha)\,f(\alpha\,y)\,f((1-\alpha)\,y)\,\dif \alpha\,.
\qedhere
\end{align*}
\end{proof}

\begin{remark}
  \label{rk:coupling-for-lemma-3.4}
  Note that the last proof makes use of a probabilistic coupling argument, since we actually prove and use the following property:
  the pair of random variables $(x\theta_x,x(1-\theta_x))$, where $\theta_x$ is distributed as $q(x,\alpha)d\alpha$, is
  stochastically increasing w.r.t.\ $x$. This means that, for all $x\leq y$, there exists a coupling of the random variables
  $\theta_x$ and $\theta_y$, i.e.\ two random variables $\theta'_x$ and $\theta'_y$ with the same laws as $\theta_x$ and $\theta_y$
  can be constructed on the same probability space, such that $x\theta'_x\leq y\theta'_y$ and $x(1-\theta'_x)\leq y(1-\theta'_y)$.
  Therefore, Assumption~\ref{hyp.fct.rep.q} means that the offspring masses of two individuals reproducing at respective masses $x$
  and $y$ can be coupled so that the masses of the offspring are in the same order as those of the parents.
\end{remark}

\medskip

\begin{proposition}
\label{da.prop.proba.extinct.decroissante}
Under Assumption \ref{hyp.fct.rep.q}, if the division rate $\taudiv(S,.)$ is non decreasing then the extinction probability $p^S:x\mapsto p^S(x)$ is non increasing.
\end{proposition}

\medskip

The assumption that $b$ increases with the mass $x$ of an individual is biologically natural, since a bigger
  total biomass usually means a bigger fraction of biomass devoted to the bio-molecular mechanisms involved in cellular division. We
  give below an analytical proof of this proposition, but it can also be proved using probabilistic arguments, as explained in
  Remark~\ref{rk:proba-prop-3.5} below.

\medskip

\begin{proof}
We prove by induction that the function $p^S_n$ is non increasing for any $n \in \NN^*$, where $p^S_n$ is  given by \eqref{proba.extinction.n.gen}.
Let $0<x<y<M$. As $A^S_u(x)<A^S_u(y)$, for any $u \geq 0$,
\begin{align*}
  p^S_1(x)
  &=
  D\,\int_0^\infty e^{-\int_0^t \taudiv(S,A^S_u(x))\,\dif u - D\,t}\, \dif t
 \geq
  D\,\int_0^\infty e^{-\int_0^t \taudiv(S,A^S_u(y))\,\dif u - D\,t}\, \dif t
	=
  p^S_1(y) \, .
\end{align*}
Then the function $p^S_1$ is non increasing. Let $n\in \NN^*$, we assume that the function $p^S_n$ is non increasing.

We can write $p^S_{n+1}(x)$ as
\begin{align*}
p^S_{n+1}(x)
	&=
		p^S_1(x)
		+
		\PP^S_{\delta_x}\left(
						\left\{\text{extinction before the $(n+1)$-th generation} \right\}
						\cap \left\{\eta^S_{\tau}\neq 0\right\}\right) \, ,
\end{align*}
with
\begin{align*}
& \PP^S_{\delta_x}\left(\left\{\text{extinction before the $(n+1)$-th generation}\right\} 
						\cap \left\{\eta^S_{\tau}\neq 0\right\}\right)
\\
	& \qquad =
		\int_{0}^\infty 
	\taudiv(S,A^S_t(x))\,e^{-\int_0^t \taudiv(S,A^S_u(x))\,\dif u-D\,t}\,
\\
	& \qquad \qquad \qquad
	\int_0^1 q(A^S_t(x),\alpha)\,p^S_n(\alpha\,A^S_t(x))\,p^S_n((1-\alpha)\,
	   A^S_t(x))\,\dif\alpha\, \dif t\,.
\end{align*}
The following relation holds
\begin{align*}
p^S_{n+1}(x)
	&=
		p^S_1(x)
	+
		 p^S_{n+1}(x \,| \,\eta^S_{\tau}\neq 0)\,(1-p^S_1(x))
\end{align*}
with
\begin{align*}
p^S_{n+1}(x \, | \,  \eta^S_{\tau}\neq 0)
	& =
	\PP^S_{\delta_x}(\text{extinction before the $(n+1)$-th generation}\, |\, \eta^S_{\tau}\neq 0)\,.
\end{align*}

\medskip
 
Since for any $t\geq 0$, $A^S_{t(x,y)+t}(x)=A^S_t(y)$, then, by a change of variable, 
\begin{align}
\nonumber
&
\int_{t_S(x,y)}^\infty 
	\taudiv(S,A^S_t(x))\,e^{-\int_0^t \taudiv(S,A^S_u(x))\,\dif u-D\,t}\,	
\\
\nonumber
	& \qquad \qquad \qquad
\int_0^1 q(A^S_t(x),\alpha)\,p^S_n(\alpha\,A^S_t(x))\,p^S_n((1-\alpha)\,
	   A^S_t(x))\,\dif\alpha\, \dif t
\\
\nonumber
	& \qquad \qquad
	=
		e^{-\int_0^{t_S(x,y)} \taudiv(S,A^S_u(x))\,\dif u-D\,t_S(x,y)}\,
		\int_{0}^\infty 
		\taudiv(S,A^S_t(y))\,e^{-\int_0^t \taudiv(S,A^S_u(y))\,\dif u-D\,t}\,
\\
\nonumber
      &\qquad\qquad\qquad\qquad\qquad\qquad\qquad\qquad
		\int_0^1 q(A^S_t(y),\alpha)\,p^S_n(\alpha\,A^S_t(y))\,
				p^S_n((1-\alpha)\,A^S_t(y))\,\dif\alpha
		\, \dif t
\\
\label{eq.I1}
	& \qquad \qquad
	=
		e^{-\int_0^{t_S(x,y)} \taudiv(S,A^S_u(x))\,\dif u-D\,t_S(x,y)} \,
		p^S_{n+1}(y \,|\, \eta^S_\tau \neq 0)\, (1-p^S_1(y)) \, .
\end{align}

For any $t\in[0,t_S(x,y)]$ we have $A^S_t(x) \leq y$. Since we assume that the function $p^S_n$ is non increasing, from Lemma \ref{lemma.int.q.dec}, we then get
\begin{align}
\nonumber
&
\int_0^{t_S(x,y)}
	\taudiv(S,A^S_t(x))\,e^{-\int_0^t \taudiv(S,A^S_u(x))\,\dif u-D\,t}\,
\\
\nonumber
	& \qquad \qquad \qquad
		\int_0^1 q(A^S_t(x),\alpha)\,p^S_n(\alpha\,A^S_t(x))\,p^S_n((1-\alpha)\,A^S_t(x))\,\dif\alpha
	\, \dif t
\\
\nonumber
	& \qquad \qquad
	\geq 
	\int_0^{t_S(x,y)}
		\taudiv(S,A^S_t(x))\,e^{-\int_0^t \taudiv(S,A^S_u(x))\,\dif u-D\,t}\,\dif t
	\int_0^1 
		q(y,\alpha)\,p^S_n(\alpha\,y)\,p^S_n((1-\alpha)\,y)\,\dif\alpha
\\
\nonumber
	& \qquad \qquad
	=
	\left(
		1-e^{-\int_0^{t_S(x,y)} \taudiv(S,A^S_u(x))\,\dif u-D\,t(x,y)}
		-D\,
		\int_0^{t_S(x,y)}
			e^{-\int_0^t \taudiv(S,A^S_u(x))\,\dif u-D\,t}\,\dif t
	\right)
\\
\label{eq.I2}
	& \qquad \qquad \qquad
	\times 
	\int_0^1 
		q(y,\alpha)\,p^S_n(\alpha\,y)\,p^S_n((1-\alpha)\,y)\,\dif\alpha \,.
\end{align}

Still because the function $p^S_n$ is non increasing and from Lemma \ref{lemma.int.q.dec},
\begin{multline*}
\PP^S_{\delta_y}
	\left(
		\left\{\text{extinction before the $(n+1)$-th generation}\right\}
		\cap 
		\left\{\eta^S_{\tau}\neq 0\right\}
	\right)
\\
	\leq
		\int_0^{\infty} 
			\taudiv(S,A^S_t(y))\,e^{-\int_0^t \taudiv(S,A^S_u(y))\,\dif u-D\,t}\, \dif t
		\int_0^1 q(y,\alpha)\,p^S_n(\alpha\,y)\,p^S_n((1-\alpha)\,y)\,\dif\alpha
\\
	=
		(1-p^S_1(y))\,
		\int_0^1 q(y,\alpha)\,p^S_n(\alpha\,y)\,p^S_n((1-\alpha)\,y)\,\dif\alpha\,.
\end{multline*}
Hence
\begin{align*}
\int_0^1 q(y,\alpha)\,p^S_n(\alpha\,y)\,p^S_n((1-\alpha)\,y)\,\dif\alpha
	& \geq p^S_{n+1}(y \,|\, \eta^S_\tau \neq 0) \,.
\end{align*}
Adding \eqref{eq.I1} and \eqref{eq.I2}, and using the last inequality, we then get
\begin{align*}
p^S_{n+1}(x \,|\, \eta^S_\tau \neq 0)
	&\geq
		\Bigg[
		\frac{1-D\,
				\int_0^{t_S(x,y)}
					e^{-\int_0^t \taudiv(S,A^S_u(x))\,\dif u-D\,t}\,
					\dif t}
			{1-p^S_1(x)}
\\
	& \qquad
		- e^{-\int_0^{t_S(x,y)} \taudiv(S,A^S_u(x))\,\dif u-D\,t_S(x,y)}\,
			\frac{p^S_1(y)}{1-p^S_1(x)}
		\Bigg]\, 
		p^S_{n+1}(y | \eta^S_\tau \neq 0)\,.
\end{align*} 
Moreover,
\begin{align*}
p^S_1(x)
	=
		D\, \int_0^{t_S(x,y)} 
				e^{-\int_0^t \taudiv(S,A^S_u(x))\,\dif u-D\,t}\,
				\dif t
		+ e^{-\int_0^{t_S(x,y)} \taudiv(S,A^S_u(x))\,\dif u-D\,t_S(x,y)} \, p^S_1(y)\,.
\end{align*}
Hence,
\begin{align*}
p^S_{n+1}(x \,|\, \eta^S_\tau \neq 0)
	& \geq
		p^S_{n+1}(y \,|\, \eta^S_\tau \neq 0)\,.
\end{align*}
Thus,
\begin{align*}
p^S_{n+1}(x) - p^S_{n+1}(y)
	& =
		p^S_1(x)+p^S_{n+1}(x \,|\, \eta^S_{\tau}\neq 0)\,(1-p^S_1(x))
\\
	& \quad
		-p^S_1(y)-p^S_{n+1}(y \,|\, \eta^S_{\tau}\neq 0)\,(1-p^S_1(y))
\\
	&\geq 
		(p^S_1(x)-p^S_1(y)) \, (1-p^S_{n+1}(y \,|\, \eta^S_{\tau}\neq 0))
	\geq  
		0 \, .
\end{align*}
This ends the induction.
Passing to the limit, we finally get
\begin{align*}
p^S(x)-p^S(y)
	&=
		\lim_{n\to \infty}
		(p^S_n(x)-p^S_n(y))
	\geq 0 \, .
\qedhere
\end{align*}
\end{proof}

\medskip

\begin{remark}
  \label{rk:proba-prop-3.5}
    The last result can also be proved by a probabilistic coupling argument as follows. First, for all $x\in(0,M)$, the time of death
    or division of an individual of mass $x$ can be constructed from an exponential random variable $E$ with parameter 1 as
    $T_x=\inf\{t\geq 0:\int_0^t (b(S,A^S_s(x))+D)ds\geq E\}$. Hence, if $x\leq y$ then $A^S_{T_x}(x)\leq A^S_{T_y}(y)$.
    Second, we observe that the probability of death given death or division occurs for an individual of mass $x$, $D/(D+b(S,x))$, is
    non-increasing as a function of $x$. Hence, using Remark~\ref{rk:coupling-for-lemma-3.4}, given $x\leq y$, we can construct a
    coupling between the branching processes $(\eta_t^S,t\geq 0)$ with $\eta^S_0=\delta_x$ and $(\hat{\eta}_t^S,t\geq 0)$ with
    $\hat{\eta}^S_0=\delta_y$ such that the random sets $M_1$ and $\hat{M}_1$ of masses at birth of the individuals of the first generation
    satisfy the following property: the cardinals $|M_1|$ and $|\hat{M}_1|$ of $M_1$ and $\hat{M}_1$ is either 0 or 2, $|\hat{M}_1|=0$
    implies that $|M_1|=0$ and if both have cardinal 2, then $M_1=\{x_1,x_2\}$ and $\hat{M}_1=\{\hat{x}_1,\hat{x}_2\}$ with
    $x_1\leq\hat{x}_1$ and $x_2\leq\hat{x}_2$.

    It then follows by induction that the processes $(\eta_t^S,t\geq 0)$ and $(\hat{\eta}_t^S,t\geq 0)$ can be coupled so that, for
    all $n\geq 0$, the masses at birth of all the individuals of the $n$-th generation can be ordered into two vectors
    $V^n=(x^n_1,\ldots,x^n_{G_n})$ and $\hat{V}^n=(\hat{x}^n_1,\ldots,\hat{x}^n_{\hat{G}_n})$, where $G_n$ and $\hat{G}_n$ are the
    random sizes of generation $n$ in $\eta^S$ and $\hat{\eta}^S$ respectively, satisfying the following property: for all $n$,
    $G_n\leq\hat{G}_n$ and for all $1\leq i\leq G_n$, $x^n_i\leq\hat{x}^n_i$. This implies
    Proposition~\ref{da.prop.proba.extinct.decroissante} since survival of $\eta^S$ means that $G_n\geq 1$ for all $n$ and this
    implies that $\hat{\eta}^S$ also survives. Hence $p^S(x)\geq p^S(y)$. 
\end{remark}

\medskip

We now extend the notation of the extinction probability with a dependence in $D$ : let $p^{S,D}(x)$ be the extinction probability of the population evolving in the environment determined by $S$, with a death rate $D$ and a initial individual with mass $x$.

\medskip

\begin{proposition}
For any $x \in [0,\mmax]$, the function $D\mapsto p^{S,D}(x)$ is non-decreasing.
\end{proposition}

\medskip

\begin{proof}
Let $D'>D$.
\begin{align*}
p_1^{S,D}(x)
	=
		D\,\int_0^\infty e^{-\int_0^t b(A^S_u(x))\,\dif u- D\,t}
	&=
		1-\int_0^\infty b(A^S_u(x))\,e^{-\int_0^t b(A^S_u(x))\,\dif u- D\,t}
\\
	& \leq
		1-\int_0^\infty b(A^S_u(x))\,e^{-\int_0^t b(A^S_u(x))\,\dif u- D'\,t}
	= 
		p_1^{S,D'}(x)\,.
\end{align*}
Hence $D\mapsto p_1^{S,D}(x)$ is non-decreasing. For $n\in\NN^*$, let assume that $D\mapsto p_n^{S,D}(x)$ is non-decreasing, then
\begin{align*}
p_{n+1}^{S,D}(x)
	&=
		1-
		\int_0^\infty b(A^S_u(x))\,e^{-\int_0^t b(A^S_u(x))\,\dif u- D\,t}\,
\\
	& \qquad \qquad \qquad
			\left[
				1-
				\int_0^1 q(A^S_t(x),\alpha)\,p_{n}^{D,S}(\alpha\,A^S_t(x))\,
					p_{n}^{D,S}((1-\alpha)\,A^S_t(x))\,\dif \alpha 
			\right]\, \dif t
\\
	&\leq
		1-
		\int_0^\infty b(A^S_u(x))\,e^{-\int_0^t b(A^S_u(x))\,\dif u- D'\,t}\,
\\
	& \qquad \qquad \qquad
			\left[
				1-
				\int_0^1 q(A^S_t(x),\alpha)\,p_{n}^{D',S}(\alpha\,A^S_t(x))\,
					p_{n}^{D',S}((1-\alpha)\,A^S_t(x))\,\dif \alpha 
			\right]\, \dif t
\\
	&=
	p_{n+1}^{S,D'}(x)
\end{align*}
Then for any $n$, $p_{n}^{S,D}(x)\leq p_{n}^{S,D'}(x)$. Passing to the limit,
\begin{align*}
	p^{S,D}(x) &= \lim_{n\to \infty} p_n^{S,D}(x) 
	\leq \lim_{n\to \infty} p_n^{S,D'}(x)=p^{S,D'}(x)\,.
	\qedhere
\end{align*}
\end{proof}

\subsection{Monotonicity properties w.r.t. $S$ on the stochastic model}
\label{subsec.mon.sto.S}

We now study the variations of the survival probability w.r.t. the environmental parameter $S$. We need additional assumptions.

\medskip

\begin{hypotheses}
\label{hypotheses.variations.FI}
\begin{enumerate}
\item \label{hyp.b.croissant} The division rate function $\taudiv$ is non decreasing in the two variables $S$ and $x$.

\item \label{hyp.g.croissant} The growth speed $g$ in non decreasing in $S$:
$$
	g(S^1,x) \leq g(S^2,x) \, , \quad \forall x \in [0,\mmax], \, 0<S^1<S^2\, .
$$

\item \label{hyp.bg.dec} For any $x \in (0,M)$, the function $S \mapsto \frac{b(S,x)}{g(S,x)}$ is non increasing.

\end{enumerate}
\end{hypotheses}

\medskip

Assumptions \ref{hyp.b.croissant} and \ref{hyp.g.croissant} above are natural from the biological point of view since a bigger total biomass means a
  bigger fraction of biomass devoted to division and a larger amount of resources means a more efficient growth and division of
  cells. 
Assumption \ref{hyp.bg.dec} means that the growth rate increases faster in $S$ that the division rate. This excludes that, increasing $S$, a faster division produces too small individuals to grow and reproduce.
  Note that these assumptions are satisfied for instance if $b$ does not depends on the variable $S$ and if $g$ is of the
form $g(S,x)=\mu(S)\,\tilde g(x)$, where $\mu$ is an non decreasing function, for example a Monod kinetics \citep{monod1949a}
$\mu(S)=\mu_{\max}\,\frac{S}{K+S}$ where $\mu_{\max}$ and $K$ are constants. The form $g(S,x)=\mu(S)\,\tilde
  g(x)$ means that the resource concentration $S$ influences the speed of growth of bacteria independently of the way $x$ influences
  growth. In other words, the flow $t\mapsto A^{S'}_t(x)$ is just a proportional time change of $t\mapsto A^S_t(x)$ for all $S,S'$.

\medskip

\begin{theorem}
\label{thm.proba.extinct.decroissante.S}
Under Assumptions \ref{hyp.fct.rep.q} and \ref{hypotheses.variations.FI}, we have for any $x\in(0,\mmax)$
$$
	\PP^{S_1}_{\delta_x}(\text{survival}) \leq \PP^{S_2}_{\delta_x}(\text{survival}) \, , \quad \forall \, 0< S^1 \leq S^2\,.
$$
\end{theorem}

\medskip

In other words, for the chemostat model, under the assumptions of the previous theorem, the higher the substrate concentration in the chemostat at the mutation time is, the higher the survival probability is.

\medskip

\begin{remark}
  \label{rk:main-probabilistic-proof}
    Following Remark~\ref{rk:proba-prop-3.5}, the last result could be also proved by probabilistic coupling arguments. These
    arguments would actually only require to assume that, for all $x\leq y$ and $S_1\leq S_2$, there exists a coupling between the
    sets $M^{x,S_1}_1$ and $M^{y,S_2}_1$ of biomasses at birth of the individuals of the first generation in the branching processes
    $(\eta^{S_1}_t,t\geq 0)$ such that $\eta^{S_1}_0=\delta_x$ and $(\eta^{S_2}_t,t\geq 0)$ such that $\eta^{S_2}_0=\delta_y$
    respectively, satisfying the following property: $|M^{y,S_2}_1|=0$ implies $|M^{x,S_1}_1|=0$ and when both have cardinal 2, then
    $M^{x,S_1}_1=\{x_1,x_2\}$ and $M^{y,S_2}_1=\{y_1,y_2\}$ with $x_1\leq y_1$ and $x_2\leq y_2$. The proof of
    Theorem~\ref{thm.proba.extinct.decroissante.S} given below actually consists in checking that the coupling assumption above is
    implied by Assumptions~\ref{hyp.fct.rep.q} and~\ref{hypotheses.variations.FI}. However, this coupling assumption is hard to check
    in practice and this is why we chose to give an analytical proof based on the Assumptions~\ref{hyp.fct.rep.q}
    and~\ref{hypotheses.variations.FI}, which are stronger, but easier to check.

    Of course, Theorem~\ref{thm.proba.extinct.decroissante.S} is certainly valid under weaker assumptions, for example if $b$ or $g$
    are not monotonic w.r.t.\ $S$, but our probabilistic approach requires coupling assumptions like the one stated in this remark,
    so the method would not extend easily to such cases.
\end{remark}



\medskip

\begin{proof}
For any $y \in (0,M)$ the function $S \mapsto g(S,y)$ is non decreasing then $A^{S^1}_u(x) \leq A^{S^2}_u(x)$ for any $u \geq 0$. Moreover the function $(S,x) \mapsto \taudiv(S,x)$ is non decreasing in the two variables $S$ and $x$, then we have
\begin{align*}
  p_1^{S^1}(x)-p_1^{S^2}(x)
	&=
		D \int_0^\infty
			e^{-D\,t}\,
			\left[
				e^{-\int_0^t \taudiv(S^1, A^{S^1}_u(x))\,\dif u}
				-
				e^{-\int_0^t \taudiv(S^2, A^{S^2}_u(x))\,\dif u}
			\right] \, \dif t
 \geq 0\,.
\end{align*}
The function $S\mapsto p_1^S(x)$ is then non increasing for any $x\in(0,\mmax)$. Let $n\in \NN^*$, we assume that the function $S\mapsto p_n^S(x)$ is non increasing for any $x\in (0,\mmax)$.

The function $t\mapsto \int_0^t (b(S,A_u^S(x))+D)\,\dif u$ is a bijection from $[0,\infty[$ to $[0,\infty[$. Hence, for $X\geq 0$, there exists a unique $T^{S}_x(X)$ such that $X=\int_0^{T^{S}_x(X)} (b(S,A_u^S(x))+D)\,\dif u$.
By the change of variable $X=\int_0^t (b(S,A_u^S(x))+D)\,\dif u$ in \eqref{proba.extinction.n.gen}, we obtain
\begin{align*}
p_{n+1}^S(x)
	&=
		\int_0^\infty \Bigg[
			\frac{D}{b\Big(S,A^{S}_{T^{S}_x(X)}(x)\Big)+D}
			+\frac{b\Big(S,A^{S}_{T^{S}_x(X)}(x)\Big)}{b\Big(S,A^{S}_{T^{S}_x(X)}(x)\Big)+D}\,
			 \Psi^{S,x}_n(X)
		\Bigg]\,e^{-X}\,\dif X
\end{align*}
with
\begin{align}
\label{def.Psi}
\Psi^{S,x}_n(X) = 
	\int_0^1 
	q\Big(A^{S}_{T^{S}_x(X)}(x),\alpha \Big)\,
	p^S_n\Big(\alpha \, A^{S}_{T^{S}_x(X)}(x) \Big)\,
	p^S_n\Big((1-\alpha)\,A^{S}_{T^{S}_x(X)}(x)\Big)\,\dif\alpha \,.
\end{align}
Moreover, for all $X\geq 0$, for $S_1\leq S_2$, by the changes of variable $A^{S^i}_u(x)=y$ for $i=1,2$ and by Assumption~\ref{hypotheses.variations.FI}-\ref{hyp.bg.dec}, we have
\begin{align*}
\int_x^{A^{S^2}_{T^{S^2}_x(X)}(x)} \frac{b(S^2,y)+D}{g(S^2,y)}\,\dif y
	= X &= \int_x^{A^{S^1}_{T^{S^1}_x(X)}(x)} \frac{b(S^1,y)+D}{g(S^1,y)}\,\dif y
\\
	&\geq 
	\int_x^{A^{S^1}_{T^{S^1}_x(X)}(x)} \frac{b(S^2,y)+D}{g(S^2,y)}\,\dif y
\end{align*}
therefore $A^{S^2}_{T^{S^2}_x(X)}(x) \geq A^{S^1}_{T^{S^1}_x(X)}(x)$. We deduce from Lemma \ref{lemma.int.q.dec} and Proposition \ref{da.prop.proba.extinct.decroissante} that 
$\Psi^{S^1,x}_n(X) \geq \Psi^{S^2,x}_n(X)$. Hence, using $\Psi^{S^1,x}_n(X) \geq \Psi^{S^2,x}_n(X)$ in the expression of $p^{S^1}_{n+1}$, subtracting the expression of $p^{S^2}_{n+1}$ and factorizing the terms, we obtain
\begin{multline*}
p_{n+1}^{S^1}(x) - p_{n+1}^{S^2}(x)
	\geq
	\int_0^\infty \Bigg[
			\frac{D}{b\Big(S^1,A^{S^1}_{T^{S^1}_x(X)}(x)\Big)+D}
			-\frac{D}{b\Big({S^2},A^{S^2}_{T^{S^2}_x(X)}(x)\Big)+D}\Bigg]
\\
		(1-\Psi^{S^2,x}_n(X))
		\,e^{-X}\,\dif X
\end{multline*}
and as $\Psi^{S^2,x}_n(X)\leq 1$, by Assumptions ~\ref{hypotheses.variations.FI}-\ref{hyp.b.croissant}, $p_{n+1}^{S^1}(x) \geq p_{n+1}^{S^2}(x)$.
Finally, passing to the limit, we get
\begin{align*}
p^{S^1}(x)-p^{S^2}(x)
	&=
		\lim_{n\to \infty}
		(p_n^{S^1}(x)-p_n^{S^2}(x))
	\geq 0 \, .
	\qedhere
\end{align*}
\end{proof}

\subsection{Properties on the variations of the eigenvalue}
\label{subsec.mon.eigen}

Until now, we only studied the probability of survival of the branching process $\eta^S$. The underlying coupling
  arguments require to consider the population state at each generation in a process where generations actually overlap. This is why
  such an approach is hard to apply directly to the integro-differential eigenvalue problem, where the notion of generations is
  difficult to define. However, the link between the stochastic and deterministic problems stated in Theorem~\ref{th.critere} allows
  to extend the monotonicity properties of Theorem~\ref{thm.proba.extinct.decroissante.S} to the eigenvalue $\Lambda_S$, as proved
  below. The next corollary is a direct consequence of Theorems \ref{th.critere} and \ref{thm.proba.extinct.decroissante.S}.
  
\medskip
  
\begin{corollary}
\label{cor.lambda.S}
Under Assumptions \ref{hyp.fct.rep.q} and \ref{hypotheses.variations.FI},
\begin{enumerate}
\item if there exists $S_1>0$ such that $\Lambda_{S_1}>0$, then $\Lambda_{S_2}>0$ for any $S_2>S_1$;
\item if there exists $S_1>0$ such that $\Lambda_{S_1}\leq 0$, then $\Lambda_{S_2}\leq 0$ for any $S_2<S_1$.
\end{enumerate}
\end{corollary} 

\medskip

This Corollary allows to deduce the following result about variation of the eigenvalue with respect to $S$.
\begin{corollary}
Under Assumptions \ref{hypotheses.variations.FI}, the function $S \mapsto \Lambda_S$ is non decreasing.
\end{corollary} 

\medskip 
The monotonicity of $b$ is important to obtain the monotonicity of the eigenvalue (and of the survival probability).
  For example, one can imagine cases where a fast growth rate $g$ transports individuals to big masses and if the division rate is
  low for high values, the monotonicity of the eigenvalue does not hold (see \cite{calvez2012a} for non monotonic examples).

\begin{proof}
Let $S^*>0$ be fixed. We set $D'=D+\Lambda_{S^*}>0$. 
Let $\Lambda_S'$ be the eigenvalue of the following eigenproblem:
$$
	\partial_x(g(S,x)\,\hat u_S'(x))+(b(S,x)+D'+\Lambda'_S)\hat u'_S(x)
	= 2\int_x^\mmax \frac{\taudiv(S,z)}{z}\, 
		q\left(z,\frac{x}{z}\right)\, 
		\hat u'_S(z)\,\dif z.
$$
For $S=S^*$, we have $\Lambda'_{S^*}=0$, then from Corollary \ref{cor.lambda.S}, for any $S\leq S^*$, $\Lambda'_{S}\leq 0$. Moreover
$$
	\Lambda'_S = \Lambda_S+D-D' =\Lambda_S-\Lambda_{S^*} \,.			
$$
Hence 
$
	\Lambda_S \leq \Lambda_{S^*} \,.
$
\end{proof}

\subsection{Extensions and concluding remarks}
\label{subsec.ext.ccl}

The previous method can be applied for more general $g$, for which the growth in one environment is larger than the growth in the other one for all masses. A particular case is given in the following corollary.

\medskip

\begin{corollary}
We assume that the division rate function $\taudiv$ does not depend on the variable $S$ and is non decreasing in the variable $x$ and that the growth speed $g$ is of the form
$
	g(S,x) = \mu(S)\,\tilde g(x)\, ,
$
where $g(S, x)>0$ for any $x \in (0,M)$ and $\tilde g \in C[0,\mmax]\cap C^1(0,\mmax)$ is such that $\tilde g(0)=\tilde g(\mmax)=0$.
Then, we have
$$\PP^{S^1}_{\delta_x}(\text{survival}) \leq \PP^{S^2}_{\delta_x}(\text{survival}) 
	\Longleftrightarrow \mu(S^1) \leq \mu(S^2)
$$
and
$$\Lambda_{S^1} \leq \Lambda_{S^2} \Longleftrightarrow \mu(S^1) \leq \mu(S^2)\,.$$
\end{corollary} 

\medskip
More generally, the following result states the link between the comparison of the survival probability and the comparison of the eigenvalue.
 
We extend the notations of the survival probability $\PP_{\delta_x}^{S, D}(\text{survival})$ and the eigenvalue $\Lambda^D_{S}$ with a dependence to the death rate $D$. 

\medskip

\begin{proposition}
\label{prop.general.method}
Let $S_1,S_2>0$. If for any $x\in[0,\mmax]$ and for any $D>0$, we have 
$\PP_{\delta_x}^{S_1, D}(\text{survival}) \geq \PP_{\delta_x}^{S_2, D}(\text{survival})$, then
$$
	\forall D>0\,, \qquad \Lambda^D_{S_1}\geq \Lambda^D_{S_2}\,.
$$
\end{proposition}

The condition on the survival probability stated in the previous theorem could be of course obtained under the
  appropriate coupling assumptions (as in Remark~\ref{rk:proba-prop-3.5}), but it seems hard to find general practical conditions on
  the parameters of the model ensuring such a property.
Note also that this coupling method could be applied for example to the case where the division distribution $q$ also depends on the
variable $S$. The results of Section \ref{subsec.mon.sto} would remain true as, for this section, the substrate concentration is fixed. The
difficulties are in the control of the variation in $S$ of $\Psi_n^{S,x}$ defined by \eqref{def.Psi}.

\medskip

\begin{proof}
Let $S_1>0$. We set $D'=D+\Lambda^D_{S_1}>0$. Let $\Lambda_S'$ be the eigenvalue associated to the eigenproblem
$$
	\partial_x(g(S,x)\,\hat u_S'(x))+(b(S,x)+D'+\Lambda'_S)\hat u'_S(x)
	= 2\int_x^\mmax \frac{\taudiv(S,z)}{z}\, 
		q\left(\frac{x}{z}\right)\, 
		\hat u'_S(z)\,\dif z\,.
$$
For $S=S_1$, $\Lambda'_{S_1}=0$, we then deduce, from Theorem \ref{th.critere}, that 
$\PP_{\delta_x}^{S_1, D'}(\text{survival})=0$ and then, by assumption, that 
$\PP_{\delta_x}^{S_2, D'}(\text{survival})=0$.
From Theorem \ref{th.critere} $\Lambda'_{S_2}\leq 0$. Moreover
$$
	\Lambda'_{S_2} = \Lambda^D_{S_2}+D-D' =\Lambda^D_{S_2}-\Lambda^D_{S_1} 			
$$
hence 
\begin{align*}
	\Lambda^D_{S_2} & \leq \Lambda^D_{S_1} \,. 
\qedhere
\end{align*}
\end{proof}


\end{document}